\documentclass{article}
\usepackage[utf8]{inputenc}
\usepackage{amsmath}
\usepackage{amssymb}
\usepackage{amsthm}

\usepackage{biblatex}

\title{Group-Graph Reciprocal Pairs}
\author{Kirsty Campbell}
\date{}

\begin{document}

\maketitle

\begin{abstract}
    In a 2018 paper, Cameron and Semeraro posed the problem of finding all group-graph reciprocal pairs. In this paper, we make a significant contribution to finding all such pairs. A group and graph form a reciprocal pair if they satisfy the relation $$P_{\Gamma,G}(x)=(-1)^nF_G(-x)$$ where $P_{\Gamma,G}(x)$ is the orbital chromatic polynomial of a graph $\Gamma$ and $F_G(x)$ is the cycle polynomial of a finite permutation group. We define a set of graphs to be \textit{$k$-stars} and prove that they satisfy a reciprocality relation with some group depending on $k$. These graphs are comprised of a complete graph with $k$ vertices and a further $\alpha$ `points' which are only connected to each vertex in the centre. This group is a subgroup of $S_k\times S_\alpha$, which is the automorphism group of a \textit{$k$-star} and $\alpha$ is the number of points on the star. We conjecture a list of group-graph reciprocal pairs. 
\end{abstract}

\section{Introduction}

The \textit{cycle polynomial} of a finite permutation group $G$ acting on a set $\Omega$ of size $n$ is $$F_G(x)=\sum_{g\in G}x^{c(g)}$$ where $c(g)$ is the number of cycles (including fixed points) of $g$ on $\Omega$. This polynomial will always be monic with degree $n$ as the only element with $n$ cycles is the identity element. We also need to define the \textit{orbital chromatic polynomial} of a graph $\Gamma$. Firstly, the \textit{chromatic polynomial} of a graph $\Gamma$, denoted by $P_{\Gamma}(x)$, is a monic polynomial of degree equal to the number of vertices of $\Gamma$, such that when evaluated at a positive integer $c$, it gives the number of proper colourings of the graph $\Gamma$ with $c$ colours. A proper colouring of a graph is a map from the vertex set to the set of colourings such that no two adjacent vertices have the same colour. \\
To get the orbital chromatic polynomial of $\Gamma$, we let $G$ be an automorphism group of $\Gamma$ and take any $g\in G$. Then define the graph denoted by $\Gamma/g$ to be the graph obtained by `contracting' each cycle of $g$ to a single vertex. The chromatic polynomial $P_{\Gamma/g}(x)$ is the polynomial which counts proper colourings of $\Gamma$ fixed by $g$. From the definition of proper colouring, it is clear to see that if any cycle of $g$ contains an edge, then $P_{\Gamma,g}(x)=0$ as there is now a loop in the graph. Now we can define the orbital chromatic polynomial of a pair $(\Gamma,G)$. It is given by the summation $$P_{\Gamma,G}(x)=\sum_{g\in G}P_{\Gamma/g}(x).$$
Note that if a graph contains loops, there are no proper colourings of this graph so we assume that all the graphs we are working with are simple.\\
We call a group-graph pair $(\Gamma,G)$ a \textit{reciprocal pair} if the relation $$P_{\Gamma,G}(x)=(-1)^nF_G(-x)$$ holds for a simple graph $\Gamma$ with $n$ vertices and where $G$ is an automorphism group of $\Gamma$. This reciprocity relationship is defined in \cite{Jason}.\\
There are already some examples of group-graph reciprocal pairs that we know of. Richard Stanley in his 1974 paper \cite{Stanley} introduced the idea of reciprocity for combinatorial polynomials. Peter Cameron and Jason Semeraro then posed the problem of finding all group-graph reciprocal pairs and we expand on their search in this paper.\\
There are of course trivial cases; the case where $\Gamma$ is a complete graph with $n$ vertices and the automorphism group $S_n$, and the null graph with any $G$ which has no odd permutations. The 4-cycle with the dihedral group of order 8 also form a group-graph reciprocal pair and the proof of this is given in \cite[Example 1]{Jason}.\\
In \cite[Theorem 4.1]{Jason}, an infinite family of group-graph reciprocal pairs is given where the graph is a star. We define a \textit{$k$-star} (with $k<n$) to be the graph with vertex set $V(\Gamma)=\left\{1,...,n\right\}$ and edge set $$E(\Gamma)=\left\{(i,j) \,|\, i\ne j, \,i,j\leq k\right\}\cup\left\{(i,j) \,|\, 1\leq i\leq k,\, k+1\leq j\leq n\right\}$$ so a 1-\textit{star} is a star. Here is the main result of our paper: 

\medskip

\noindent\textbf{Theorem 1} \textit{Suppose $(\Gamma, G)$ is a reciprocal pair where $\Gamma$ is a $k$-star with $n\geq 2k+1$ vertices and $G<S_n$. Then the following hold: 
\begin{enumerate}
    \item $n=r(k+1)+k$ for some $r\in\mathbb N$ 
    \item $G=\bar{G}\times S_k$ for some $\bar{G}\leq S_{r(k+1)}$ where $n=r(k+1)+k$ and either 
    \begin{enumerate}
        \item $(S_{k+1})^r\leq\bar{G}\leq S_{k+1}\wr S_r$ if $k$ is odd or
        \item $(S_{k+1})^r\leq\bar{G}\leq S_{k+1}\wr A_r$ if $k$ is even.
    \end{enumerate}
\end{enumerate}}

\medskip

\noindent \cite[Theorem 4.1]{Jason} speaks to the case $k=1$. Its proof is also similar, however in the proof of our result we are using methods such as counting arguments giving the proof in this paper some novelty.\\
This gives us a doubly infinite family of reciprocal pairs. It is also known that we can form more pairs by taking direct products and wreath products and was proved in \cite[Proposition 3.5, 3.6]{Jason}. For example, if we have $n$ reciprocal pairs $(\Gamma_1,G_1),(\Gamma_2,G_2),...,(\Gamma_n,G_n)$ and let $\Gamma$ be the disjoint union of $\Gamma_1,...,\Gamma_n$ and $G$ be the direct product of $G_1,...,G_n$ where $G_i\leq \Gamma_i$. Then $(\Gamma,G)$ is also a reciprocal pair. A similar thing happens when we take the wreath product. Let $\Gamma$ be the graph given by the disjoint union of $m$ copies of a graph $\Gamma_0$ and $G\leq \Gamma_0$. Assume that $(\Gamma_0,G)$ form a reciprocal pair. If we have another group $H$ of permutations of degree $m$ containing no odd permutations, then $(\Gamma,G\wr H)$ is also a reciprocal pair. This means that we can form new pairs from already existing ones.

\medskip

\noindent\textbf{Definition:}\textit{A reciprocal pair $(\Gamma,G)$ is \emph{irreducible} is there is no non-trivial decomposition $\Gamma=\Gamma_1\cup\Gamma_2\cup...\cup\Gamma_n$ and subgroups $G_i\leq G$ with $G_i\leq\text{Aut}(\Gamma_i)$ such that either \begin{enumerate}
    \item $G=G_1\times...\times G_n$ with each $(\Gamma_i,G_i)$ a reciprocal pair,
    \item $G_i\cong G_j$ and $\Gamma_i\cong \Gamma_j$ for all $i\ne j$ and $G=G_0\wr H$ for some $H\leq S_n$ and $G_0\cong G_i$ for all $i$.
\end{enumerate}}

\noindent To find these pairs, we conducted a search using the software Magma \cite{Magma} to speed up the process. During this search, the only pairs that we found we either (products of) \textit{$k$-stars}, the 4-cycle and dihedral group of order 8, and the two trivial pairs. This leads us to the following conjecture:

\medskip

\noindent\textbf{Conjecture:} \textit{Let $(\Gamma,G)$ be an irreducible reciprocal pair. Then one of the following hold: \begin{enumerate}
    \item $(\Gamma,G)$ is one of the two trivial cases: \begin{enumerate}
        \item $\Gamma$ is a null graph and $G$ is any group with no odd permutations,
        \item $\Gamma$ is the complete graph with $n$ vertices and $G=S_n$.
    \end{enumerate}
    \item $\Gamma$ is a 4-cycle and $G\cong \text{Aut}(\Gamma)$.
    \item $\Gamma$ is a $k$-star and $G$ is defined in Theorem 1.
\end{enumerate}}

\section{Proof of Theorem 1}

To prove Theorem 1, we first state some results from \cite{Jason} which are used in the proof. 

\medskip

\noindent\textbf{Lemma 2} \cite[Proposition 1.3]{Jason} \textit{If $G$ contains no odd permutations, then $F_G$ is an even or odd function according as $n$ is even or odd; in other words, $$F_G(-x)=(-1)^nF_G(x).$$}

\medskip

\noindent\textbf{Lemma 3} \cite[Proposition 1.6]{Jason} \textit{Suppose that $G$ and $H$ are permutation groups on the same set with $H\leq G$. Suppose that $F_H(-a)=0$, for some positive integer $a$. Then also $F_G(-a)=0$. }

\medskip

\noindent\textbf{Lemma 4} \cite[Proposition 1.7]{Jason} \textit{Let $G_1$ and $G_2$ be permutation groups on disjoint sets $\Omega_1$ and $\Omega_2$ respectively. Let $G=G_1\times G_2$ acting on $\Omega_1\cup\Omega_2$ . Then
$$F_G(x)=F_{G_1}(x)\cdot F_{G_2}(x).$$}

\medskip

\noindent\textbf{Lemma 5} \cite[Proposition 1.8]{Jason} \textit{$$F_{G\wr H}(x)=|G|^mF_H\left(F_G(x)/|G|\right),$$ where $m$ is the size of the domain of $H$.}

\medskip

\noindent\textbf{Lemma 6} \cite[Lemma 3.1]{Jason} \textit{Suppose that $(\Gamma,G)$ is a reciprocal pair. Then the number of edges of $\Gamma$ is the sum of the transpositions in $G$ and the number of transpositions $(i,j)\in G$ for which $\left\{i,j\right\}$ is a non-edge. In other words we have $E(\Gamma)=t_(G)+t_0(G)$ where $t_0(G)$ denotes the transpositions which are non-edges.}

\medskip

\noindent\textbf{Proposition 7} \textit{Let $\Gamma$ be a $k$-star as defined above. Then the full automorphism group $G$ of $\Gamma$ is isomorphic to $S_k\times S_\alpha$ where $\alpha$ is the number of points on the $k$-star and $S_\alpha$ acts on the points around the star.}

\medskip

\noindent\textbf{Proof} Clearly, any element $g\in S_k$ or $g\in S_\alpha$ is also an element of any automorphism group of $\Gamma$. It remains to show that there does not exist an element $g\in G$ such that it contains a cycle containing a vertex from the centre of the star and vertex which is a point of the star. By definition, $g$ is an element of the automorphism group if and only if $\left\{i,j\right\}$ and $\left\{g(i),g(j)\right\}$ are both edges in $\Gamma$. However if $g$ contains a cycle satisfying the above, then the vertex which is a point of the star gets mapped to the centre and is now connected to other (non-centre) points of the star. This gives a contradiction, as there are no edges between points of $\Gamma$ and therefore $g$ cannot be an element of any automorphism group of $\Gamma$.

\medskip

\noindent Now we have the necessary background lemmas, we are in a position to prove Theorem 1. 

\medskip

\noindent\textbf{Lemma 8} \textit{$(\Gamma,G)$ is a reciprocal pair if and only if \begin{equation}
    F_{\bar{G}}(-x)F_T(-x)=(-1)^n x(x-1)\cdots(x-k+1)F_{\bar{G}}(x-k)
\end{equation}
where $\bar{G}=stab_G(1,...,k)$, $T\leq S_k$ and $G=\bar{G}\times T$.}

\medskip

\noindent\textit{Proof} From the definition of group-graph reciprocal pairs, $(\Gamma,G)$ is a reciprocal pair if and only if $F_G(-x)=(-1)^nP_{\Gamma,G}(x)$. The chromatic polynomial of a \textit{$k$-star} is $$P_{\Gamma,G}(x)=\sum_{g\in G}P(\Gamma/g)(x)=\sum_{g\in \bar{G}\times S_1^k}x(x-1)\cdots(x-k+1)(x-k)^{c(g)-k}.$$ In order to avoid any elements that gives $P(\Gamma/g)=0$ we must restrict the sum to any $g\in\bar{G}\times S_1^k$. This gives $$P_{\Gamma,G}(x)=\frac{x(x-1)\cdots(x-k+1)}{(x-k)^k}\sum_{g\in \bar{G}\times S_1^k}(x-k)^{c(g)}$$ $$=\frac{x(x-1)\cdots(x-k+1)}{(x-k)^k}F_{\bar{G}\times S_1^k}(x-k).$$ By Lemma 4 we get $$F_{\bar{G}\times S_1^k}(x-k)=F_{\bar{G}}(x-k)F_{S_1^k}(x-k)=F_{\bar{G}}(x-k)(x-k)^k.$$ Hence
$$P_{\Gamma,G}(x)=x(x-1)\cdots(x-k+1)F_{\bar{G}}(x-k).$$ Then, again using Lemma 4 with $G=\bar{G}\times T$ we get $$F_{\bar{G}}(-x)F_T(-x)=(-1)^n x(x-1)\cdots(x-k+1)F_{\bar{G}}(x-k).$$

\medskip

\noindent\textbf{Lemma 9} \textit{$G=\bar{G}\times S_k$ for some $\bar{G}$ such that $(S_{k+1})^r\leq \bar{G}\leq S_{k+1}\wr H$ for some $H\leq S_r$ where $n=r(k+1)+k$.}

\medskip

\noindent\textbf{Proof} Let $\alpha$ be the number of points on the \textit{$k$-star}. Then the number of edges is $k\alpha+\frac{1}{2}k(k-1)$. By Lemma 6, $t(G)+t_0(G)=k\alpha+\frac{1}{2}k(k-1)$ where $t_0(G)=t(G)-\delta$ for some $0\leq\delta\leq\frac{1}{2}k(k-1)$. One of the consequences of Proposition 7 is that there can be no element in any automorphism group $G$ of $\Gamma$ which acts between the centre of the star and its points. Hence, the only place a transposition which is also an edge can exist is in the centre of the star which has $\frac{1}{2}k(k-1)$ edges, giving the equality for $\delta$. \\
We can also get another equation for $t_0(G)$ by taking the group generated by all the elements in $t_0(G)$: $$\langle g\in t_0(G)\rangle\cong S_{\lambda_1}\times S_{\lambda_2}\times\cdots\times S_{\lambda_r}$$ where $\sum^r_{i=1}\lambda_i=\alpha$. This is because the only place in the \textit{$k$-star} where there are no edges between vertices are between the points of the star. Also note that this group lies within $\bar{G}$. So now we have two equations for $t_0(G)$ to work with. 
$$t_0(G)=\frac{1}{2}\left(k\sum^r_{i=1}\lambda_i+\frac{1}{2}k(k-1)-\delta\right)$$
$$t_0(G)=\frac{1}{2}\left(\sum^r_{i=1}\lambda_i^2-\sum^r_{i=1}\lambda_i\right)$$
Combining these into one equation yields: 
\begin{equation}
   \Rightarrow (k+1)\sum^r_{i=1}\lambda_i+\frac{1}{2}k(k-1)-\delta=\sum^r_{i=
   1}\lambda_i^2 
\end{equation}\\
Now let there be some $\lambda_i\geq k+2$ for some $0\leq i\leq r$. Then this gives $S_{k+2}\leq\bar{G}$ and hence $F_{\bar{G}}(x)$ will have a root at $-(k+1)$. By Lemma 8 with $x=-1$ we have $$F_{\bar{G}}(1)F_T(1)=(-1)^n(-1)(-2)\cdots(-k)F_{\bar{G}}(-1-k)$$
$$\Rightarrow |\bar{G}||T|=(-1)^n(-1)(-2)\cdots(-k)\cdot 0$$
which is a contradiction, so $\lambda_i\leq k+1$ for all $i$.\\
Now define $p_v:=|\left\{\lambda_v : \lambda_v=v\right\}|$. Note that $p_i\geq0$ for all $i$. This makes equation (2) read as $$(k+1)p_1+\cdots+ k(k+1)p_k+(k+1)^2p_{k+1}+\frac{1}{2}k(k-1)-\delta $$$$=p_1+2^2p_2+\cdots+(k+1)^2p_{k+1}$$\\
Moving all the terms to the left hand side gives $$kp_1+\cdots+i(k+1-i)p_i+\cdots+kp_k+\frac{1}{2}k(k-1)-\delta=0$$\\
As all $p_i\geq0$ and $\frac{1}{2}k(k-1)-\delta\geq0$ we must have that all $p_i=0$ for all $i\leq k$. This gives us a few conclusions. First, that $\delta=\frac{1}{2}k(k-1)$ so we must have that $S_k\leq T$. But before we had that $T\leq S_k$ so $T=S_k$ and $G=S_k\times\bar{G}$. Second, we have $p_{k+1}=r$ to satisfy the condition $\sum^r_{i=1}\lambda_i=\alpha$, so $\alpha=r(k+1)$ and $n=r(k+1)+k$ which proves the first part of Theorem 1. Third, we have that $(S_{k+1})^r\leq\bar{G}\leq S_{k+1}\wr H$ for some $H\leq S_r$ and the lemma is proved.  

\medskip

\noindent\textbf{Proof of Theorem 1} The previous lemma proved the first part of Theorem 1. It also proved that $T=S_k$, so now we can substitute $F_T(x)=x(x+1)+...+(x+k-1)$ into equation (1) and get $$F_{\bar{G}}(-x)(-1)^kx(x-1)\cdots(x-k+1)=(-1)^nx(x-1)\cdots(x-k+1)F_{\bar{G}}(x-k)$$ $$\Rightarrow F_{\bar{G}}(-x)=(-1)^{r(k+1)}F_{\bar{G}}(x-k)$$ which holds if and only if $(\Gamma,G)$ is a reciprocal pair. \\
As seen in Lemma 5, the cycle polynomial of $\bar{G}$ can also be written in terms of the wreath product $$F_{\bar{G}}=((k+1)!)^rF_H\left(\frac{F_{S_{k+1}}(x)}{(k+1)!}\right).$$ Given that we know the cycle polynomial of $S_{k+1}(x)=x(x+1)\cdots(x+k)$ we also have $S_{k+1}(-x)=(-1)^{k+1}x(x-1)\cdots(x-k)$ and $S_{k+1}(x-k)=x(x-1)\cdots(x-k)$. Using those cycle polynomials, we can use Lemma 5 to find $F_{\bar{G}}(-x)$ and $F_{\bar{G}}(x-k)$: $$F_{\bar{G}}(-x)=((k+1)!)^rF_H\left(\frac{(-1)^{k+1}x(x-1)\cdots(x-k)}{(k+1)!}\right)$$
$$=((k+1)!)^rF_H\left((-1)^{k+1}y\right)$$
$$F_{\bar{G}}(x-k)=((k+1)!)^rF_H\left(\frac{x(x-1)\cdots(x-k)}{(k+1)!}\right)$$
$$=(k+1)!)^rF_H(y)$$\\
Now let $H\leq S_r$ and let $k$ be odd (part 2(a) of Theorem 1) and observe that $F_{\bar{G}}(-x)=F_H(y)=F_{\bar{G}}(x-k)$. Hence, as $k+1$ is even, we have that $(\Gamma,G)$ is a reciprocal pair. \\
For the last part of Theorem 1, let $H\leq A_r$ and let $k$ be even. Then we get $$F_{\bar{G}}(-x)=F_H(-y)=(-1)^{r(k+1)}F_H(y)=(-1)^{r(k+1)}F_{\bar{G}}(x-k).$$ By Lemma 2 this holds if $H$ contains no odd permutations, but as $H\leq A_r$ that is clearly the case, so again we have that $(\Gamma,G)$ is a reciprocal pair and the theorem is proved. 
\begin{flushright}
$\Box$
\end{flushright}

\noindent This gives a new infinite family of group-graph reciprocal pairs. As discussed before, we know we can form new pairs by taking wreath and direct products of existing pairs. As of yet, we have not found any pairs not in the form described in the conjecture in the first section. Future research into this area could investigate into the conjecture and try and find a proof that all group-graph reciprocal pairs are of the form described. 

\medskip

\medskip

\begin{center}
    \textbf{Acknowledgements}
\end{center}
 I would like to thank the London Mathematical Society for funding this project through the Undergraduate Research Bursary as well as my supervisor Dr. Jason Semeraro for suggesting the project. I've enjoyed working with him and I am extremely grateful for his support and guidance throughout this process. I would also like to thank Dr. Frank Neumann for his helpful comments and feedback on this paper. Lastly, I would like to thank Prof. Jeremy Levesley for putting me in touch with the opportunity and his ongoing support throughout this project.

\printbibliography

@article{Jason,
author={Peter J. Cameron and Jason Semeraro},
title={The Cycle Polynomial of a Permutation Group},
journal={Electronic J. Combinatorics},
volume={25},
number={1},
year={2018}
}

@article{Stanley,
author={Richard Stanley},
title={Combinatorial Reciprocity Theorems},
journal={Combinatorics},
year={1974},
pages={194-253}
}

@article {Magma,
    AUTHOR = {Bosma, Wieb and Cannon, John and Playoust, Catherine},
     TITLE = {The {M}agma algebra system. {I}. {T}he user language},
      NOTE = {Computational algebra and number theory (London, 1993)},
   JOURNAL = {J. Symbolic Comput.},
  FJOURNAL = {Journal of Symbolic Computation},
    VOLUME = {24},
      YEAR = {1997},
    NUMBER = {3-4},
     PAGES = {235--265},
      ISSN = {0747-7171},
   MRCLASS = {68Q40},
  MRNUMBER = {MR1484478},
       DOI = {10.1006/jsco.1996.0125},
       URL = {http://dx.doi.org/10.1006/jsco.1996.0125},
}
\end{document}